\documentclass[11pt]{article}
\usepackage{amsfonts,amsmath,latexsym,color,epsfig}
\date{}

	\newcommand{\RR}{{\mathbb R}}
	\newcommand{\PP}{{\mathbb P}}

\title{Disjointness graphs of segments in the space}

\author{
{\sl J\'anos Pach}\thanks{
R\'enyi Institute, Budapest and MIPT, Moscow; \texttt{pach@renyi.hu}; \texttt{pach@cims.nyu.edu}.
Supported by the National Research, Development and Innovation Office (NKFIH) project KKP-133864, ERC Advanced Grant ``GeoScape,'' the Austrian Science Fund grant Z 342-N31, and by the Ministry of Education and Science of the Russian Federation in the framework of MegaGrant No.\ 075-15-2019-1926.}
\and
{\sl G\'abor Tardos}\thanks{
R\'enyi Institute, Budapest  and Central European University, Budapest; \texttt{tardos@renyi.hu}.
Supported by the Cryptography ``Lend\"ulet'' project of the Hungarian Academy
of Sciences and by the National Research, Development and Innovation Office (NKFIH)
projects K-116769, K-132696, and KKP-133864, by the ERC Synergy Grant ``Dynasnet'' No.\ 810115, by the Ministry of Education and Science of the Russian Federation in the framework of MegaGrant No.\ 075-15-2019-1926 and by ERC Advanced Grant ``GeoScape,''.}
\and
{\sl G\'eza T\'oth}\thanks{
R\'enyi Institute, Budapest; \texttt{geza@renyi.hu}.
Supported by National Research, Development and Innovation Office, NKFIH, K-131529, KKP-133864, NKFIH-1158-6/2019, and by the Higher Educational Institutional Excellence Program 2019 and by ERC Advanced Grant ``GeoScape,''.}}

\date{}

\begin{document}

\maketitle

\begin{abstract}
The {\em disjointness graph} $G=G({\cal S})$
of a set of segments ${\cal S}$ in
${\RR}^d$, $d\ge 2,$ is a graph whose vertex set is
${\cal S}$ and two vertices are connected by
an edge if and only if the corresponding segments are disjoint. We prove that
the chromatic number of $G$ satisfies $\chi(G)\le(\omega(G))^4+(\omega(G))^3$,
where $\omega(G)$ denotes the clique number of $G$. It follows that $\cal S$ has
$\Omega(n^{1/5})$ pairwise intersecting or pairwise disjoint
elements. Stronger bounds are established for lines in space, instead of
segments.

We show that computing $\omega(G)$ and $\chi(G)$ for disjointness graphs of
lines in space are NP-hard tasks. However, we can design efficient
algorithms to compute proper colorings of $G$ in which the number of colors
satisfies the above upper bounds. One cannot expect similar results for sets
of continuous arcs, instead of segments, even in the plane. We construct
families of arcs whose disjointness graphs are triangle-free ($\omega(G)=2$),
but whose chromatic numbers are arbitrarily large.
\end{abstract}

\section{Introduction}

Given a set of (geometric) objects, their {\em intersection graph} is a graph
whose vertices correspond to the objects, two vertices being connected by an
edge if and only if their intersection is nonempty. Intersection graphs of
intervals on a line~\cite{H57}, more generally, chordal graphs~\cite{B61,D61},
and comparability graphs~\cite{D50}, turned out to be {\em perfect graphs},
that is, for them and for each of their induced subgraph $H$, we have
$\chi(H)=\omega(H)$, where $\chi(H)$ and $\omega(H)$ denote the chromatic
number and the clique number of $H$, respectively. It was shown~\cite{HS58}
that the complements of these graphs are also perfect. Based on the above
results, Berge~\cite{B61} conjectured and Lov\'asz~\cite{L72} proved that the
complement of every perfect graph is perfect.
%By now, we have a complete
%characterization of all perfect graphs~\cite{CRST06}, which immediately
%implies the Lov\'asz theorem.
\smallskip

Most geometrically defined intersection graphs are not perfect. However, in
many cases they still have nice coloring properties. For example, Asplund and
Gr\"unbaum \cite{AG60} proved that every intersection graph $G$ of
axis-parallel rectangles in the plane satisfies
$\chi(G)=O((\omega(G))^2)$. It is not known if the stronger
bound $\chi(G)=O(\omega(G))$ also holds for these graphs.
For intersection graphs of chords of a circle,
Gy\'arf\'as~\cite{G85} established the bound
$\chi(G)=O((\omega(G))^24^{\omega(G)})$, which was improved to
$O(2^{\omega(G)})$ in~\cite{KK97} and, recently, to $O(\omega\log\omega)$ by Chalermsook and Walczak \cite{ChW19}. Here we have examples of graphs $G$ with $\chi(G)$ slightly
superlinear in $\omega(G)$ \cite{K88}.
In some cases, there is no functional dependence between
$\chi$ and $\omega$. The first such example was found by Burling~\cite{B65}:
there are sets of axis-parallel boxes in ${\RR}^3$, whose intersection graphs
are {\em triangle-free} ($\omega=2$), but their chromatic numbers are
arbitrarily large. Following Gy\'arf\'as and Lehel~\cite{GL83}, we call a
family $\cal G$ of graphs {\em $\chi$-bounded} if there exists a function $f$
such that all elements $G\in\cal G$ satisfy the inequality $\chi(G)\le
f(\omega(G))$. The function $f$ is called a {\em bounding function} for $\cal
G$. Heuristically, if a family of graphs is $\chi$-bounded, then its members
can be regarded ``nearly perfect". Consult \cite{GL85,G87,K04} for
surveys.
\smallskip

At first glance, one might believe that, in analogy to perfect graphs, a
family of intersection graphs is $\chi$-bounded if and only if the family of
their complements is. Burling's above mentioned constructions show that this
is not the case: the family of complements of intersection graphs of
axis-parallel boxes in ${\RR}^d$ is $\chi$-bounded with bounding function
$f(x)=O(x\log^{d-1}x)$, see~\cite{K91}.
More recently, Pawlik, Kozik, Krawczyk, Laso\'n, Micek,
Trotter, and Walczak~\cite{PKK14} have proved that Burling's triangle-free
graphs can be realized as intersection graphs of segments in the
plane. Consequently, the family of these graphs is {\em not}
$\chi$-bounded either. On the other hand, the family of their complements
is, see Theorem~0.
\smallskip

To simplify the exposition, we call the complement of the intersection graph
of a set of objects their {\em disjointness graph}. That is, in the
disjointness graph two vertices are connected by an edge if and only if the
corresponding objects are disjoint. Using this terminology, the following
is a direct consequence of a result of
Larman, Matou\v sek, Pach, and T\"or\H ocsik.

\smallskip

\noindent{\bf Theorem 0.} \cite{LMPT94} {\em The family of disjointness graphs
  of segments in the plane is $\chi$-bounded. More precisely, every such graph
  $G$ satisfies the
inequality $\chi(G)\le(\omega(G))^4$.}

\smallskip

For the proof of Theorem~0, one has to introduce {four} partial orders on the
family of segments, and apply Dilworth's theorem~\cite{D50} four
times. Actually, the same argument also works systems of arbitrary $x$-monotone curves,
instead of segments. (A continuous curve is called $x$-monotone if any line parallel to the $y$-axis
intersects it in at most one point.) It was proved in \cite{PaT19} that in this setting the order
of magnitude or the upper bound $(\omega(G))^4$ cannot be improved.

One of the main results of this paper is a generalization of Theorem 0 to higher dimensions.
We establish
the following.
\medskip

\noindent {\bf Theorem 1.}
{\em The disjointness graph $G$ of any system of segments in ${\RR}^d, d\ge 2$
satisfies the inequality $\chi(G)\le(\omega(G))^4+(\omega(G))^3$.

Moreover, there is a polynomial time algorithm that, given the segments
corresponding to the vertices of $G$, finds a
complete subgraph $K\subseteq G$ and a proper coloring of $G$ with at most
$|V(K)|^4+|V(K)|^3$ colors.}

\medskip
Unfortunately, the technique applied in the plane does not seem to work in
higher dimensions.

If we consider full lines in place of segments, we obtain stronger bounds.

\medskip

\noindent{\bf Theorem 2.}
{\em (i) Let $G$ be the disjointness graph of a set of lines in ${\RR}^d, \; d\ge 3.$
Then we have $\;\; \chi(G)\le(\omega(G))^3.$
\smallskip

(ii) Let $G$ be the disjointness graph of a set of lines in the projective
space ${\PP}^d, \; d\ge 3.$ Then we have $\;\;  \chi(G)\le(\omega(G))^2.$
\smallskip

In both cases, there are polynomial time algorithms that, given the lines
corresponding to the vertices of $G$, find
complete subgraphs $K\subseteq G$ and proper colorings of $G$ with at most $|V(K)|^3$ and
$|V(K)|^2$ colors, respectively.}
\medskip

Note that the difference between the two scenarios comes from the fact that
parallel lines in the Euclidean space are disjoint, but the corresponding lines
in the projective space intersect.

\smallskip

Answering a question in an earlier draft of this paper, Norin \cite{N17} showed that the family of {\em intersection graphs} of lines in ${\RR}^d$ or ${\PP}^d$, $d\ge 3$, is {\em not} $\chi$-bounded.

\smallskip

Most computational problems for geometric intersection and disjointness graphs
are hard. It was shown by Kratochv\'\i l and Ne\v set\v ril~\cite{KN90} and
by Cabello, Cardinal, and Langerman~\cite{CCL13} that finding the clique
number $\omega(G)$ resp.\ the independence number $\alpha(G)$ of disjointness
graphs of segments in the plane are NP-hard. It is also known that
computing the chromatic number $\chi(G)$ of
disjointness and intersection graphs of segments in the plane is
NP-hard~\cite{EET86}. Our next theorem shows that some of the analogous problems are also
NP-hard for disjointness graphs of lines in space, while others are
tractable in this case. In particular, according to Theorem 3(i), in a
disjointness graph $G$ of lines, it is NP-hard to determine
$\omega(G)$ and $\chi(G)$. In view of this, it is interesting that one
can design polynomial time algorithms to find proper colorings and
complete subgraphs in $G$, where the number of
colors is bounded in terms of the size of the complete subgraphs,
in the way specified in the closing statements of Theorems~1 and 2.

\medskip

\noindent{\bf Theorem 3.} {\em
(i) Computing the clique number $\omega(G)$ and the chromatic number $\chi(G)$
of disjointness graphs of lines in $\RR^3$ or in $\PP^3$ are NP-hard
problems.
\smallskip

(ii) Computing the independence number $\alpha(G)$ of disjointness graphs of
lines in $\RR^3$ or in $\PP^3$, and deciding for a fixed $k$ whether
$\chi(G)\le k$, can be done in polynomial time.}

%(iii) To decide the 3-colorability of intersection graphs of lines in $\RR^3$ is an NP-complete problem.}
\medskip

The bounding functions in Theorems~0, 1, and~2 are not likely to be
optimal. As for Theorem~2 (i), we will prove that there are
disjointness graphs $G$ of lines in $\mathbb{R}^3$ for which $\frac{\chi(G)}{\omega(G)}$
are arbitrarily large. Our best constructions for disjointness graphs $G'$ of lines in the projective space
satisfy $\chi(G')\ge2\omega(G')-1$; see Theorem~2.3.
\smallskip

The proof of Theorem~1 is based on Theorem~0. Any strengthening of Theorem~0
leads to improvements of our results. For example,
if $\chi(G)=O((\omega(G))^\gamma)$ holds with any $3\le \gamma\le 4$ for the
disjointness graph of every set of segments in the plane, then the proof of
Theorem~1 implies the same bound for disjointness graphs of segments in higher
dimensions. In fact, it is sufficient to verify this statement in $3$
dimensions. For $d\ge 4$, we can find a projection in a generic direction
to the $3$-dimensional
space that does not create additional intersections and then we can apply the
$3$-dimensional bound. We focus on the case
$d=3$.
\smallskip

It follows immediately from Theorem 0 that the disjointness (and, hence,
the intersection) graph of any system of $n$ segments in the plane has a
clique or an independent set of size at least $n^{1/5}$. Indeed, denoting by
$\alpha(G)$ the maximum number of independent vertices in $G$,
we have
$$\alpha(G)\ge \frac{n}{\chi(G)}\ge\frac{n}{(\omega(G))^4},$$
so that $\alpha(G)(\omega(G))^4\ge n$. Analogously, Theorem~1 implies that
$\max(\alpha(G),\omega(G))\ge (1-o(1))n^{1/5}$ holds for disjointness (and
intersection) graphs of segments in any dimension $d\ge 2$. For disjointness
graphs of $n$ lines in $\mathbb{R}^d$ (respectively, in $\mathbb{P}^d$), we obtain
that $\max(\alpha(G),\omega(G))$ is $\Omega(n^{1/4})$ (resp.,
$\Omega(n^{1/3})$).
%In particular, the latter statement implies the bound $\Omega(n^{1/3})$
%even for lines in $\mathbb{R}^d$, provided that no two of
%them are parallel.
Using more advanced algebraic techniques, Cardinal, Payne,
and Solomon~\cite{CPS16} proved the stronger bounds
$\Omega(n^{1/3})$ (resp.,
$\Omega(n^{1/2})$).
%the same bound without
%any assumption.
\smallskip

If the order of magnitude of the bounding functions in Theorems~0 and~1 are improved,
then the improvement carries over to the lower bound on
$\max(\alpha(G),\omega(G))$. Despite many
efforts~\cite{LMPT94,KPT97,K12} to construct intersection graphs of planar
segments with small clique and independence numbers, the best known construction,
due to Kyn\v cl \cite{K12}, gives only
$$\max(\alpha(G),\omega(G))\le n^{\log 8/\log 169}\approx n^{0.405},$$
where $n$ is the number of vertices. This bound is roughly the square of the
best known lower bound.

Our next theorem shows that any improvement of the lower bound on
$\max(\alpha(G),\omega(G))$ in the plane, even if it was not achieved by an improvement of
the bounding function in Theorem~0, would also carry over to higher dimensions.

\medskip

\noindent{\bf Theorem~4.}
{\em If the disjointness graph of any set of $n$ segments in the plane has a
  clique or an independent set of size $\Omega(n^\beta)$ for some fixed
  $\beta\le1/4$, then the same is true for disjointness graphs of segments in
${\RR}^d$ for any $d>2$.}

\medskip

A continuous arc in the plane is called a {\em string}. One may wonder whether
Theorem~0 can be extended to disjointness graphs of strings in place of
segments.
The answer is no, in a very strong sense.

\medskip

\noindent{\bf Theorem~5.}
{\em There exist triangle-free disjointness graphs of $n$ strings in the plane
  with arbitrarily large chromatic numbers. Moreover, we can assume that these
  strings are simple polygonal
paths consisting of at most $4$ segments.}

\medskip

Recently, M\"utze, Walczak, and Wiechert~\cite{MWW17} modified our constructions to obtain families
 of polygonal paths consisting of {\em three} segments each, which meet the requirements of Theorem 5 and any two paths cross at most once. It is not known whether one can find polygonal paths consisting of only {\em two} segments satisfying these conditions.

The following problems remain open.

\medskip

\noindent{\bf Problem~6.} {\em (i) Is the family of disjointness graphs
of
polygonal paths, each consisting of at most two segments,
%strings in
%the plane, any pair of which intersect in at most one point,
$\chi$-bounded?
\smallskip

(ii) Is the previous statement true under the additional assumption that any
two of the polygonal paths intersect in at most one point?}
%
%the strings are polygonal paths, each consisting of at most $k$ segments, where $k>1$ is fixed?}

%---------------------------------------------------------------------------

\smallskip

This paper is organized as follows. In the next section, we prove Theorem~2,
which is needed for the proof of Theorem~1. Theorem~1 is established in
Section~3. The proof of Theorem 4 is presented in Section 4. In
Section~5, we construct several examples of disjointness graphs whose chromatic
numbers are much larger than their clique numbers. In particular,
we prove Theorem~5 and some similar statements. The last section contains the
proof of Theorem~3 and remarks on the computational complexity of related problems.

%-----------------------------------------------------------------------------------

\section{Disjointness graphs of lines--Proof of Theorem 2}

\noindent {\bf Claim 7.} {\em Let $G$ be the disjointness graph of a set of $n$
lines in
$\PP^d$.
If $G$ has an isolated vertex,
then $G$ is perfect.}
\medskip

\noindent{\bf Proof.} Let $\ell_0\in V(G)$ be a line representing an
isolated vertex of $G$.
Consider the bipartite multigraph $H$ with vertex set $V(H)=A\cup B$, where $A$
consists of all points of $\ell_0$ that belong to at least one other line
$\ell\in V(G)$, and $B$ is the set of all ($2$-dimensional) planes passing
through $\ell_0$ that contain at least one other line $\ell\in V(G)$ different
from $\ell_0$. We associate with any line $\ell\in V(G)$ different from $\ell_0$ an edge $e_\ell$ of $H$, connecting the point $p=\ell\cap\ell_0\in A$ to
the plane $\pi\in B$ that contains $\ell$. Note that there may be several
parallel edges in $H$. See Figure 1.

%\smallskip

%\begin{figure}[ht]
%\begin{center}
%\scalebox{0.5}{\includegraphics{parosgraf.eps}}
%\caption{Construction of graph $H$ in the proof of Claim 7.}
%\end{center}
%\end{figure}

\begin{figure}[ht]
\begin{center}
\scalebox{0.45}{\includegraphics{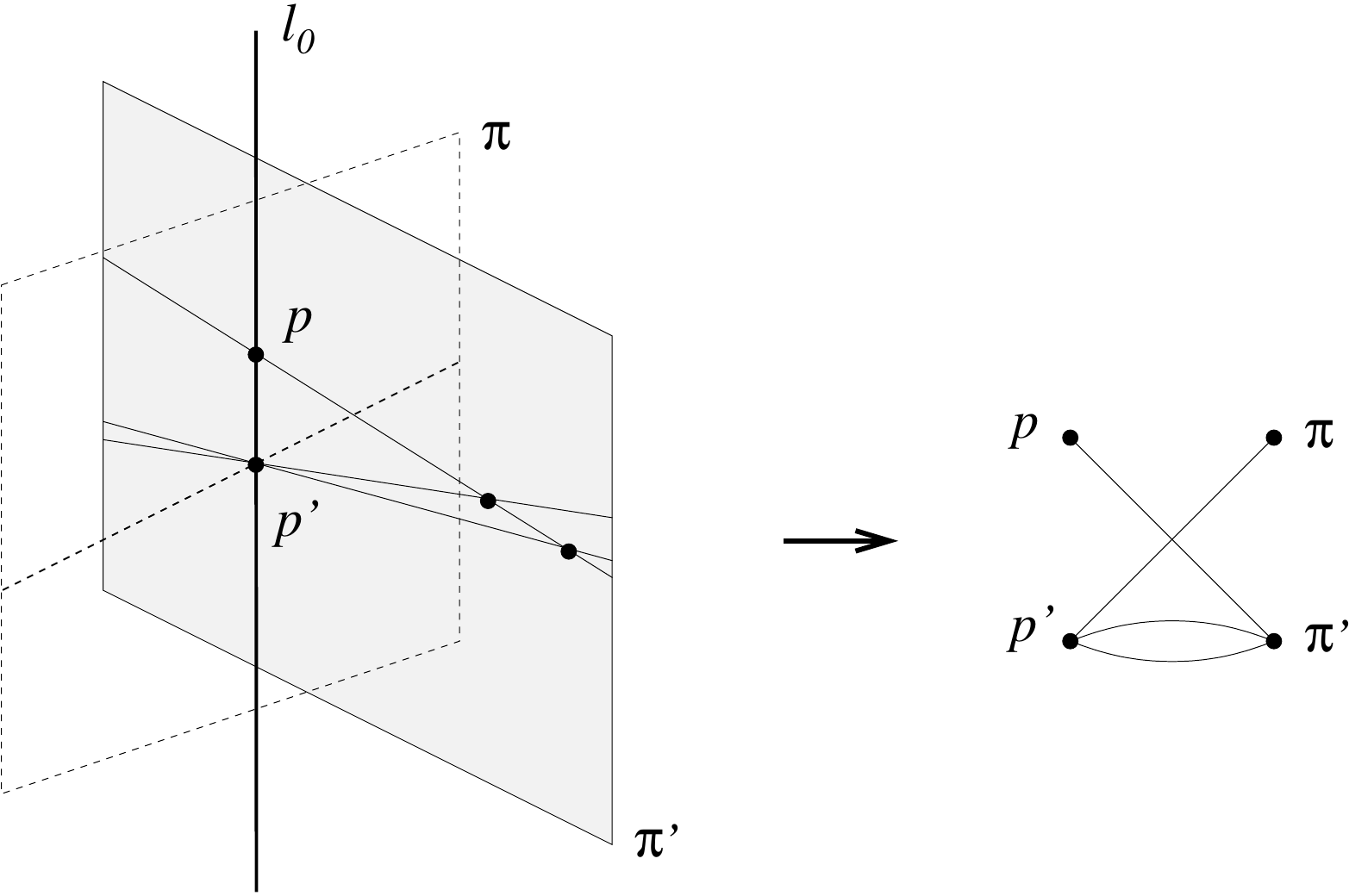}}
\caption{Construction of graph $H$ in the proof of Claim 7.}
\end{center}
\end{figure}

Observe that two lines $\ell,\ell'\in V(G)\setminus\{\ell_0\}$ intersect if
and only if $e_\ell$ and $e_{\ell'}$ share an endpoint. This means that $G$
minus the isolated vertex $\ell_0$ is isomorphic to the complement of the line
graph of $H$. The line graphs of bipartite multigraphs and their
complements are known to be perfect. (For the complements of line graphs, this
is the K\"onig-Hall theorem; see, e. g., \cite{L93}.)
The graph $G$ can be obtained by adding the isolated vertex $\ell_0$ to a perfect graph, and is,
 therefore, also
perfect. $\Box$

\medskip

\noindent{\bf Proof of Theorem 2.} We start with the proof of part (ii).  Let
$G$ be a disjointness graph of lines in $\PP^d$. Let $C\subseteq G$ be a
maximal clique
in $G$. Clearly, $|C|\le\omega(G)$. By the maximality of $C$,
for every $\ell\in V(G)\setminus C$,
there exists $c\in C$ that is not adjacent to $\ell$ in $G$. Hence,
there is a partition of $V(G)$ into disjoint sets $V_c, c\in C,$ such that
$c\in V_c$ and $c$ is an isolated vertex in the induced subgraph $G[V_c]$ of
$G$. Applying Claim 7 separately to each subgraph $G[V_c]$, we obtain
$$\chi(G)\le\sum_{c\in C}\chi(G[V_c])=\sum_{c\in C}\omega(G[V_c])
\le|C|\omega(G)\le(\omega(G))^2.$$
\smallskip

Now we turn to the proof of part (i) of Theorem 2. Let $G$ be a disjointness graph of lines
in $\RR^d$. Consider the lines in $V(G)$ as lines in the projective space
$\PP^d$, and consider the disjointness graph $G'$ of these projective
lines. Clearly, $G'$ is a subgraph of $G$ with the lines
$\ell$, $\ell'\in V(G)$ adjacent in $G$
but not adjacent in $G'$ if and only if $\ell$ and $\ell'$ are parallel. Thus,
an independent set in $G'$ induces a disjoint union of complete subgraphs in $G$, where
the vertices of each complete subgraph correspond to pairwise parallel lines.
If $k$ is the maximal number
of pairwise parallel lines in $V(G)$, then $k\le\omega(G)$ and each independent
set in $G'$ can be partitioned into at most $k$ independent sets in
$G$. Applying part (ii), we obtain
$$\chi(G)\le k\chi(G')\le \omega(G)(\omega(G'))^2\le(\omega(G))^3.$$

Finally, we prove the last claim concerning polynomial time algorithms. In the
proof of part (ii), we first took a maximal clique $C$ in $G$. Such a clique
can be efficiently found by a greedy algorithm. The partition of $V(G)$ into
subsets $V_c, c\in C,$ such that $c\in V_c$ is an isolated vertex in the
subgraph $G[V_c]$, can also be done efficiently. It remains to find a clique
of maximum size and a proper coloring of each perfect graph $G[V_c]$ with the
smallest number of colors.
It is well known that for perfect graphs, both of these tasks can be completed
in polynomial time.
See e.g. Corollary 9.4.8 on page 298 of \cite{GLS88}.
Alternatively, notice that in the proof of Claim~7 we
showed that $G[V_c]$ is, in fact, the
complement of the line graph of a bipartite multigraph (plus an isolated
vertex). Therefore, finding a maximum size complete subgraph corresponds to
finding a maximum size matching in a bipartite graph, while finding an optimal
proper coloring of $G[V_c]$ corresponds to finding a minimal size vertex cover
in a bipartite graph. This can be accomplished by much simpler and faster algorithms than
the general purpose algorithms developed for perfect graphs.

To finish the proof of the algorithmic claim for part (ii), we can simply
output as $K$ the set $C$ or one of the largest maximum cliques in $G[V_c]$
over all $c\in C$, whichever is larger. We color each $V_c$ optimally, with
pairwise disjoint sets of colors.

For the algorithmic claim about part (i), first color the
corresponding arrangement of projective lines, and then refine the coloring by
partitioning each color class into at most $k$ smaller classes, where $k$ is
the maximum number of parallel lines in the arrangement. It is easy to find
the value of $k$, just partition the lines into groups of parallel lines.
Output as $K$ the set we found for the projective lines, or a set of $k$ parallel lines,
whichever is larger. $\Box$

\medskip

\noindent{\bf Theorem 8.} {\em (i) There exist disjointness graphs $G$
of families of lines in $\mathbb{R}^3$ for which the ratio $\chi(G)/\omega(G)$ is arbitrarily large.

(ii) For any $k$ one can find a system of lines in
$\mathbb{P}^3$ whose disjointness graph $G$ satisfies $\omega(G)=k$ and $\chi(G)=2k-1$.}

\medskip

\noindent{\bf Proof.} First, we prove (i).
For some $m$ and $d$ to be determined later, consider the set $W_m^d$ of
integer points in the $d$-dimensional hypercube $[1,m]^d$.
That is, $W_m^d=\{ 1, 2, \ldots , m\}^d$.
A {\em combinatorial line} is a sequence of $m$
distinct points of $x^1,\ldots, x^m\in W_m^d$ such that for every fixed
$1\le i\le d$,
the $i$th coordinates of $x^j$, $(x^j)_i$, are either the same for all $1\le j\le m$ or
we have $(x^j)_i=j$ for all $1\le j\le m$. Note that the points of any
combinatorial line lie on a geometric straight line. Let ${\cal L}$ denote the
set of these geometric lines.

\smallskip

Let $G$ denote the disjointness graph of $\cal L$.
Since each line in $\cal L$ passes through $m$ points of $W_m^d$, and
$|W_m^d|=m^d$, we have
$\omega(G)\le m^{d-1}$. (It is easy to see that equality holds here,
but we do not need this fact for the proof.)
\smallskip

Consider any proper coloring of $G$. The color classes are families of
pairwise crossing lines in $\cal L$. Observe that any such family has a
common point in $W_m^d$, except some families consisting of  $3$ lines. Take
an optimal proper coloring of $G$ with $\chi(G)$ colors, and split each
$3$-element color class into two smaller classes. In the resulting coloring,
there are at most $2\chi(G)$ color classes, each of which has a point of $W_m^d$ in
common. This means that the set of at most $2\chi(G)$ points of $W_m^d$ (the ``centers'' of
the color classes) ``hits'' every combinatorial line. By the density version of the
Hales--Jewett theorem, due to Furstenberg and Katznelson \cite{B98, FK91}, if
$d$ is large enough relative to $m$, then any set containing fewer than half of the
points of $W_m^d$ will miss an entire combinatorial line. Choosing any
$m$ and a sufficiently large $d$ depending on $m$, we conclude that $2\chi(G)\ge
m^d/2$ and $\chi(G)/\omega(G)\ge m/4$.

\smallskip

Note that the family $\cal L$ consists of lines in $\RR^d$. To find a similar
family in 3-space, simply take the image of $\cal L$ under a projection to
$\RR^3$. One can pick a generic projection that does not change the
disjointness graph $G$. This completes the proof of part (i). Note that the
same construction does not work for projective lines, as the combinatorial
lines in $W_m^d$ fall into $2^d-1$ parallel classes, so the chromatic number of the corresponding
projective disjointness graph is smaller than $2^d$.

\smallskip

To establish part (ii), fix a positive integer $k$, and consider a set $S$ of
$2k+1$ points
in general position (no four in a plane) in $\mathbb{R}^3\subseteq
\mathbb{P}^3$. Let $\cal L$
denote the set of
$\binom{2k+1}{2}$ lines determined by them. Note that by the general
position assumption, two lines in $\cal L$ intersect if and only if
they have a point of $S$ in common. This means that the
disjointness
graph $G$ of $\cal L$ is isomorphic to the {\em Kneser graph} $G^{*}(2k+1,2)$ formed by
all $2$-element subsets of a $(2k+1)$-element set. Obviously,
$\omega(G^{*}(n,m))=\lfloor n/m\rfloor$, so $\omega(G)=k$. By a celebrated result of
Lov\'asz \cite{L78}, $\chi(G^*(n,m))=n-2m+2$  for all $n\ge2m-1$. Thus, we have $\chi(G)=2k-1$, as claimed. $\Box$

%-------------------------------------------------------------------

\section{Disjointness graphs of segments--Proof of Theorem 1}

If all segments lie in the same plane, then by Theorem~0 we
have $\chi(G)\le(\omega(G))^4$. Our next theorem generalizes this result to
the case
where the segments lie in a bounded number of distinct planes.

\medskip

\noindent{\bf Theorem 9.} {\em Let $G$ be the disjointness graph of a set of
segments in $\mathbb{R}^d, d>2,$ that lie in the union of $k$ two-dimensional
planes. We have
$$\chi(G)\le(k-1)\omega(G)+(\omega(G))^4.$$

Given the segments representing the vertices of $G$ and $k$ planes
containing them, there is a polynomial time algorithm to find a complete subgraph
$K\subseteq G$ and a proper coloring of $G$ with at most
$(k-1)|V(K)|+|V(K)|^4$ colors.}

\medskip

\noindent{\bf Proof.}  Let $\pi_1, \pi_2, \ldots, \pi_k$ be the planes
containing the segments. Partition the vertex set of $G$ into the classes
$V_1, V_2, \ldots, V_k$ by putting a segment $s$ into the class $V_i$, where
$i$ is the {\em largest}
index for which $\pi_i$ contains $s$.
\smallskip

For $i=1, 2,\ldots, k,$ we define subsets $W_i, Z_i\subseteq V_i$ with
$Z_i\subseteq W_i\subseteq V_i$ by a recursive procedure, as follows. Let
$W_1=V_1$ and let $Z_1\subseteq W_1$ be a {\em maximal size} clique in $G[W_1]$.

Assume that the sets $W_1,\dots, W_i$ and $Z_1,\ldots, Z_i$ have already been
defined for some $i<k$. Let $W_{i+1}$ denote the set of all vertices in
$V_{i+1}$ that are adjacent to every vertex in $Z_1\cup Z_2\cup\ldots\cup
Z_i$, and let $Z_{i+1}$ be a maximal size clique in
$G[W_{i+1}]$. By definition, $\bigcup_{i=1}^kZ_i$ induces a complete subgraph
in $G$, and we have
$$\sum_{i=1}^k|Z_i|\le\omega(G).$$
\smallskip

Let $s$ be a segment belonging to $Z_i$, for some $1\le i< k$. A point $p$ of
$s$ is called a {\em piercing point} if $p\in \pi_j$ for some $j>i$. Notice
that in this case, $s$ ``pierces'' the plane $\pi_j$ in a single point,
otherwise we would have $s\subset\pi_j$, contradicting our assumption that
$s\in V_i$. Letting $P$ denote the set of piercing points of all segments in
$\bigcup_{i=1}^kZ_i$, we have
$$|P|\le\sum_{i=1}^k(k-i)|Z_i|\le(k-1)\sum_{i=1}^k|Z_i|\le(k-1)\omega(G).$$
\smallskip

Let $V_0=V(G)\setminus\bigcup_{i=1}^kW_i$. We claim that every segment in
$V_0$ contains at least one piercing point. Indeed, if $s\in V_i\setminus W_i$
for some $i\le k$, then $s$ is not adjacent in $G$ to at least one segment
$t\in Z_1\cup\ldots\cup Z_{i-1}$. Thus, $s$ and $t$ are not disjoint, and
their intersection point
is a piercing point, at which $t$ pierces the plane $\pi_i$.
\smallskip

Assign a color to each piercing point $p\in P$. Coloring every segment in
$V_0$ by the color of one of its piercing points, we get a proper coloring of
$G[V_0]$ with $|P|$ colors, so that $\chi(G[V_0])\le|P|.$
\smallskip

For every $i\le k$, all segments of $W_i$ lie in the plane $\pi_i$. Therefore,
we can apply Theorem~0 to their disjointness graph $G[W_i]$, to conclude that
$\chi(G[W_i])\le(\omega(G[W_i]))^4$. By definition, $Z_i$ induces a maximum
complete subgraph in
$G[W_i]$, hence $|Z_i|=\omega(G[W_i])$ and $\chi(G[W_i])\le|Z_i|^4$.
\smallskip

Putting together the above estimates, and taking into account that
$\bigcup_{i=1}^kZ_i$ induces a complete subgraph in $G$,
we obtain
$$\chi(G)\le\chi(G[V_0])+\sum_{i=1}^k\chi(G[W_i])\le|P|+\sum_{i=1}^k|Z_i|^4$$
$$\le(k-1)\omega(G)+(\sum_{i=1}^k|Z_i|)^4\le(k-1)\omega(G)+(\omega(G))^4,$$
as required.

We can turn this estimate into a polynomial time algorithm as required,
using the fact that the proof of Theorem~0 is constructive. In particular, we use
that, given a family of segments in the plane, one can efficiently find a subfamily
$K$ of pairwise disjoint
segments and a proper coloring of the disjointness graph with at
most $|K|^4$ colors. This readily follows from the proof of Theorem~0, based on
the four easily computable (semi-algebraic) partial orders on the family of segments,
introduced in \cite{LMPT94}.

Our algorithm finds the sets $V_i$, as in the proof. However, finding $W_i$
and a maximum size clique $Z_i\subseteq W_i$ is a challenge. Instead, we use the
constructive version of Theorem~0 to find a $Z_i\subseteq W_i$ which is a
clique but not
necessarily the largest and
a proper coloring of $G[W_i]$. The definition of $W_i$ remains unchanged.
Next, the algorithm identifies the piercing points.

The algorithm outputs the clique $K=\bigcup Z_i$
and the coloring of $G$. The latter one is obtained by combining the previously
constructed colorings of the subgraphs $G[W_i]$ (using disjoint sets of colors for
different subgraphs), and coloring each remaining vertex by a previously unused
color, associated with one of the piercing points the corresponding segment passes through.
$\Box$
\smallskip

\noindent{\bf Proof of Theorem 1.} Consider the set of all lines in the {\em
  projective} space $\mathbb{P}^d$ that contain at least one segment belonging
to $V(G)$. Let $\bar{G}$ denote the disjointness graph of these
lines. Obviously, we have
$\omega(\bar{G})\le\omega(G)$. Thus, Theorem 2(ii) implies that
$$\chi(\bar{G})\le (\omega(\bar{G}))^2\le(\omega(G))^2.$$
\smallskip

Let $C$ be the set of lines corresponding to the vertices of a maximum complete subgraph in
$\bar{G}$. Fix an optimal proper coloring of $\bar{G}$.
Suppose that we used $k$  ``planar'' colors
(each such color is given to a set of lines that lie in the same plane)
and $\chi(\bar{G'})-k$ ``pointed'' colors (each given to the vertices
corresponding to a set of lines passing
through a common point).
%
%We find a
%certain number, $k$, of ``planar'' colors
%(each such color is given to a set of lines that lie in the same plane) and the remaining
%$\chi(\bar{G'})-k$ colors are ``pointed'' (each given to the vertices
%corresponding to a set of lines passing
%through a common point).
\smallskip

Consider now $G$, the disjointness graph of the segments. Let $G_0$ denote the
subgraph of $G$ induced by the set of segments whose supporting lines received
one of the $k$ planar colors in the above coloring of $\bar{G}$. These
segments lie in at most $k$ planes. Therefore, applying Theorem~9 to $G_0$,
we obtain
$$\chi(G_0)\le (k-1)\omega(G_0)+(\omega(G_0))^4
\le (k-1)\omega(G)+(\omega(G))^4.$$

For $i, 1\le i\le \chi(\bar{G})-k,$ let $G_i$ denote the
subgraph of $G$ induced by the set of segments whose supporting lines are
colored by the $i$\/th pointed color.
Its complement, $\overline{G_i}$, can be represented as the intersection graph
of subtrees of a tree. Therefore, by the result of Gavril \cite{G74},
$\overline{G_i}$ is a chordal graph, that is, it
contains no induced cycle of length
larger than $3$.
According to a theorem of Hajnal and Sur\'anyi~\cite{HS58},
any graph with this property is perfect, consequently $G_i$ is also perfect.
Therefore,
$$\chi(G_i)=\omega(G_i)\le\omega(G).$$

Putting these bounds together, we obtain that
$$\chi(G)\le \chi(G_0)+\sum_{i=1}^{\chi(\bar{G})-k}\chi(G_i)
\le(k-1)\omega(G)+(\omega(G))^4+\sum_{i=1}^{\chi(\bar{G})-k}\omega(G)$$
$$\le((\omega(\bar{G}))^2-1)\omega(G)+(\omega(G))^4<(\omega(G))^3
+(\omega(G))^4.$$
\smallskip

To prove the algorithmic claim in Theorem~1, we first apply the
algorithm of Theorem~2 to the disjointness graph $\bar G'$. We distinguish between
planar and pointed color classes and find the subgraphs $G_i$. We output
a coloring of $G$, where for each $G_i, i>0$ we use the smallest possible number of colors
($G_i$ is perfect, so its optimal coloring can be found in polynomial time), and
we color $G_0$ by the algorithm described in Theorem~9. The subgraphs $G_i$ are
colored using pairwise disjoint sets of colors. We output the largest clique
$K$ that we can find. This may belong to a subgraph $G_i$ with $i>0$, or may be
found in $G_0$ or in $\bar G'$ by the algorithms given by Theorem~9 or
Theorem~2, respectively. (In the last case, we need to turn a clique in
$\bar G'$ into a clique of the same size in $G$, by picking an
arbitrary segment from each of the pairwise disjoint lines.)
$\Box$

%------------------------------------------------------------------------------------

%---------------------------------------------------------------

\section{Ramsey-type bounds in $\mathbb{R}^2$ vs.
$\mathbb{R}^3$--Proof of Theorem 4}

As we have pointed out in the Introduction, it is sufficient to establish
Theorem 4 in $\mathbb{R}^3$. We rephrase Theorem 4 for this case in the
following form.
\medskip

\noindent{\bf Theorem 10.} {\em Let $f(m)$ be a function with the property
  that for any disjointness graph $G$ of a system of segments in
  $\mathbb{R}^2$ with
$\max(\alpha(G),\omega(G))\le m$ we have $|V(G)|\le f(m).$

Then for any disjointness graph $G$ of a system of segments in
$\mathbb{R}^3$
with
$\max(\alpha(G),\omega(G))$
$\le m$ we have $|V(G)|\le f(m)+m^4.$}

\medskip

Applying Theorem 10 with $f(k)=ck^{1/\beta}$, Theorem 4 immediately
follows. We prove
Theorem~10 by adapting the proof of Theorem~9.
\medskip

\noindent{\bf Proof of Theorem 10.} Let $G$ be the disjointness graph of a
set of
segments in $\mathbb{R}^3$ with $\omega(G)\le m$ and $\alpha(G)\le m$.
\smallskip

First, assume that all segments lie in the union of $k$ planes, for some
$k\ge 1$. Define the sets of vertices $V_i$, $W_i$, and $Z_i$ for every
$1\le i\le k$, as in the proof of Theorem~9, and let
$V_0=V(G)\setminus\bigcup_{i=1}^kW_i$.
Since all elements of $W_i$ lie in the same plane,
the subgraph induced by them is a planar segment disjointness
graph for every $i\ge1$. We can clearly represent these graphs by segments in
a common plane $\pi$ such that two segments intersect if and only they come
from the same set $W_i$ and there they intersect.
In this way, we obtain a system of segments in the
plane whose disjointness graph $G^*$ is the {\em join} of the graphs $G[W_i]$,
i.e., $G^*$ is obtained by taking the disjoint union of $G[W_i]$ (for all
$i\ge1$) and adding all edges between $W_i$ and $W_j$ for every pair $i\ne
j$.
Clearly, we have
$$\omega(G^*)=\sum_{i=1}^k\omega(G[W_i])=\sum_{i=1}^k|Z_i|\le\omega(G)\le m,$$
and $$\alpha(G^*)=\max_{i=1}^k\alpha(G[W_i])\le\alpha(G)\le m.$$

By our assumption, $G^*$ has at most $f(m)$ vertices, so that
$\sum_{i=1}^k|W_i|\le f(m).$  As we have seen in the proof of Theorem
9, the total number of piercing points is at most
$(k-1)\sum_{i=1}^k|Z_i|\le(k-1)\omega(G)<km$, and each segment in $V_0$
contains at least one of them. Each piercing point is contained in at most $m$
segments, because these segments induce an independent set in $G$.
Thus, we have $|V_0|<km^2$ and
$$|V(G)|=|V_0|+\bigcup_{i=1}^k|W_i|<km^2+|V(G^*)|\le km^2+f(m).$$
\medskip

Now we turn to the general case, where there is no bound on the number of
planes containing the segments. As in the proof of Theorem~1, we consider the
disjointness graph $\bar{G}$ of the supporting lines of the segments in the
projective space $\mathbb{P}^3$. Clearly, we have
$\omega(\bar{G})\le\omega(G)\le m$,
so by
Theorem~1 we have $\chi(\bar{G})\le m^2$. Following the proof of Theorem~1,
take an optimal coloring of $\bar{G}$, and let $G_0$ denote the subgraph of
$G$ induced by the segments whose supporting lines received one of the planar
colors. Letting $k$ denote the number of planar colors, for every $i, 1\le
i\le \omega(\bar{G})-k,$ let $G_i$ denote the subgraph of $G$ induced by the
set of segments whose supporting lines received the $i$th pointed color. As in
the proof of Theorem 1, every $G_i, i\ge 1$ is perfect and, hence, its number
of vertices satisfies
$$|V(G_i)|\le\chi(G_i)\alpha(G_i)\le\omega(G_i)\alpha(G_i)\le\omega(G)\alpha(G)\le
m^2.$$
The segments belonging to $V(G_0)$ lie in at most $k$ planes. In view of the
previous paragraph, $|V(G_0)|\le km^2+f(m)$ vertices.
Combining the above bounds, we obtain
$$|V(G)|=|V(G_0)|+\sum_{i=1}^{\chi(\bar{G})-k}|V(G_i)|
\le km^2+f(m)+(\chi(\bar{G})-k)m^2$$
$$\le km^2+f(m)+(m^2-k)m^2\le f(m)+m^4,$$
which completes the proof.   $\Box$

%-----------------------------------------------------------------------------

\section{Constructions--Proof of Theorem 5}

%Shift graphs were introduced by Erd\H os and Hajnal \cite{EH64}.
The aim of this section is to describe various arrangements of geometric
objects in 2, 3, and 4 dimensions with triangle-free disjointness graphs,
whose chromatic numbers grow logarithmically with the number of objects. (This
is much faster than the rate of growth in Theorem 8.) Our constructions can
be regarded as geometric realizations of a sequence of graphs discovered by
Erd\H os and Hajnal.

\medskip

\noindent{\bf Definition 11.} \cite{EH64}. {\em Given $m>1$, let
$H\noindent_m$, the {\em $m$th shift graph}, be a graph whose vertex set
consists of all ordered pairs $(i, j)$ with  $1\le i < j\le m$, and two
pairs $(i, j)$ and $(k, l)$ are connected by an edge if and only if $j=k$ or
$l=i$.}

\medskip

Obviously, $H_m$ is triangle-free for every $m>1$. It is not hard to show
(see, e.g., \cite{L93}, Problem 9.26)
that $\chi(H_m) = \lceil \log_2m\rceil$.
Therefore, Theorem 5 follows directly from part (vii) of the next theorem.

\medskip

\noindent {\bf Theorem 12.}
{\em For every $m$, the shift graph $H_m$ can be obtained as a disjointness
graph,
where each vertex is represented by

(i) a line minus a point in $\RR^2$;

(ii) a two-dimensional plane in $\RR^4$;

(iii) the intersection of two general position half-spaces in $\RR^3$;

(iv) the union of two segments in $\RR^2$;

(v) a triangle in $\RR^4$;

(vi) a simplex in $\RR^3$;

(vii) a polygonal curve in $\RR^2$, consisting of four line segments.}

\medskip

\def\ell{L}
\noindent{\bf Proof.}
(i) Let $\ell_1, \ldots , \ell_m$ be lines in general position in the plane.
For any $1\le i < j\le m$,  let us represent the pair $(i, j)$
by the ``pointed line'' $p_{ij}=\ell_i\setminus\ell_j$.

Fix $1\le i < j\le m$,   $1\le k < l\le m$, and set $X=p_{ij}\cap
p_{kl}=(\ell_i\cap \ell_k)\setminus(\ell_j\cup\ell_l)$. If $i=k$, then $X$ is an
infinite set.

Otherwise, $\ell_i\cap\ell_k$ consists of a single point. In this case, $X$ is
empty if and only if this point belongs to $\ell_j\cup\ell_l$. By the general
position assumption, this happens if and only if $j=k$ or $l=i$. Thus, the
disjointness graph of the sets $p_{ij},\; 1\le i<j\le m,$ is isomorphic to the
shift graph $H_m$.

\smallskip

(ii) Let $h_1, \ldots, h_m$ be hyperplanes in general position in $\RR^4$.  For
every $i$, fix another hyperplane $h'_i$, parallel (but not identical)
to $h_i$.
For any $1\le i < j\le m$,  represent the pair $(i, j)$
by the two dimensional plane $p_{ij}=h_i\cap h'_j$.

Given $1\le i < j\le m$,   $1\le k < l\le m$,
the set $X=p_{ij}\cap p_{kl}=h_i\cap h'_j\cap h_k\cap h'_l$ is the
intersection of four hyperplanes. If the four hyperplanes are in general
position, then $X$ consists of a single point.

If the hyperplanes are not in general position, then some of the four indices
must coincide. If $i=k$ or $j=l$, then two of the hyperplanes coincide and $X$
is a line. In the remaining cases, when $j=k$ or $l=i$,
among the four hyperplanes two are parallel, so their intersection $X$ is empty.

\smallskip

(iii) For $i=1,\dots,m$, define the half-space $h_i$ as
$$h_i=\{(x,y,z)\in\RR^3\mid ix+i^2y+i^3z>1\}.$$
Note that the bounding planes of these half-spaces are in
general position. For any $1\le i < j\le m$, represent the pair $(i, j)$
by $p_{ij}=h_j\setminus h_i$.

Now let $1\le i < j\le m$,   $1\le k < l\le m$. If $j=k$ or $l=i$, the
sets $p_{ij}$ and $p_{kl}$ are obviously disjoint. If $i=k$ or $j=l$, then
$p_{ij}\cap p_{kl}$ is the intersection of at most 3 half-spaces in general
position, so it is unbounded and not empty.

It remains to analyze the case when all four
indices are distinct. This requires some calculation. We assume
without loss of generality that $j<l$. Consider the point $P=(x,y,z)\in\RR^3$
with $x=\frac1i+\frac1j+\frac1k$, $y=-\frac1{ij}-\frac1{jk}-\frac1{ki}$ and
$z=\frac1{ijk}$. This is the intersection point of the bounding planes of
$h_i$, $h_j$ and $h_k$. Therefore, the polynomial $zu^3+yu^2+xu-1$ vanishes at
$u=i,j,k$, and it must be positive at
$u=l$, as $l>i,j,k$ and the leading coefficient is positive. This means that
$P$ lies in the open half-space
$h_l$. As the bounding planes of $h_i$, $h_j$ and $h_k$ are in general
position, one can find a point $P'$ arbitrarily close to $P$ (the intersection
point of these half-planes) with $P'\in
h_j\setminus(h_i\cup h_k)$. If we choose $P'$ close enough to $P$, it will also
belong to $h_l$. Thus, $P'\in p_{ij}\cap p_{kl}$, and so $p_{ij}$ and
$p_{kl}$ are not disjoint.
\smallskip

(iv), (v), and (vi) directly follow from (i), (ii) and (iii), respectively, by
replacing the unbounded geometric objects representing the vertices with their sufficiently large bounded subsets.
\smallskip

(vii) Let $C$ be an almost vertical, very short curve (arc) in the plane,
convex from the right (that is, the set of points to the right of $C$ is
convex)
lying in a small neighborhood of $(0,1)$.
%That is, it forms a convex set together with the point $(1, 0)$.
Let $p_1, p_2, \ldots , p_{m}$ be a sequence of $m$ points on $C$
such that $p_j$ is
above $p_i$ if and only if $j>i$.
For every $1\le i\le m$, let $T_i$ be an equilateral triangle whose base is
horizontal, whose upper vertex is $p_{i}$, and whose center is on the $x$-axis.
Let $q_i$ and $r_i$ be the lower right and lower left vertices of $T_i$, respectively.
It is easy to see that $T_j$ contains $T_i$ in its interior if $j>i$.
Let $s_i$ be a point on $r_ip_{i}$, very close to $p_i$.

Let us represent the vertex $(i,j)$ of the shift graph $H_m$ by the polygonal
curve
$p_{ij}=t_{ij}p_jq_jr_js_j$, where the point $t_{ij}$ is on the $x$-axis slightly
to the left of the line $p_ip_j$. Note that if $C$ is short enough and close
enough to vertical, then $t_{ij}$ can be chosen so that it belongs to the
interior of all triangles $T_k$ for $1\le k\le m$. In particular, the
entire polygonal path $p_{ij}$ belongs to $T_j$.

It depends on our earlier choices of
the vertices $p_{i'}$, how close we have to choose $s_i$ to
$p_i$. Analogously, it depends on our earlier choices of $p_{i'}$ and
$s_{i'}$, how close we have to choose $t_{ij}$ to the line. Instead of describing
an explicit construction, we simply claim that with proper choices
of these points, we obtain a disjointness representation of the shift graph.

\smallskip

To see this, let $1\le i < j\le m$,   $1\le k < l\le m$. If $j=l$, then three
of the four line segments in $p_{ij}$ and $p_{kl}$ are the same, so they
intersect. Otherwise, assume without loss of generality that $j<l$. As noted
above, $p_{ij}$ belongs to the triangle $T_j$, which, in turn, lies
in the interior of $T_l$. Three segments of $p_{kl}$ lie on the edges of $T_l$,
so if $p_{ij}$ and $p_{kl}$ meet, the fourth segment,
$t_{kl}p_l$, must meet $p_{ij}$. This segment enters the triangle $T_j$, so it meets one of
its edges. Namely, for $j>k$ it follows from the convexity of the curve $C$
that the
segment $t_{kl}p_l$ intersects the edge $p_jq_j$ and, hence, also $p_{ij}$.
Analogously, if $j<k$, then $t_{kl}p_l$ intersects the interior of the edge
$r_jp_j$. This is true even if $t_{kl}$ were chosen {\em on} the line
$p_kp_l$, so choosing $s_j$ close enough to $p_j$, one can make sure that
$t_{kl}p_l$ intersects $r_js_j$ and, hence, also $p_{ij}$.
On the other hand, if $j=k$, we choose $t_{kl}$ so that $t_{kl}p_l$ is
just slightly to the left of $p_j=p_k$, so it enters $T_j$ through the
interior of the segment $s_jp_j$ that is {\em not} contained in $p_{ij}$. To
see that in this case $p_{ij}$ and $p_{kl}$ are disjoint, it is enough to check
that $t_{kl}p_l$ and $t_{ij}p_j$ are disjoint. This is true, because $p_j$ is on
the right of $t_{kl}p_l$ and (from the convexity of $C$) the slope of the
segments is such that $p_j$ is the closest point of the segment $t_{ij}p_j$ to $t_{kl}p_l$.
$\Box$

\section{Complexity issues--Proof of Theorem 3}

The aim of this section is to outline the proof of Theorem 3 and to establish
some related complexity results. For simplicity, we only consider systems of
lines in the {\em projective} space $\PP^3$. It is easy to see that by
removing a generic hyperplane (not containing any of the intersection points),
we can turn a system of projective lines into a system of lines into $\RR^3$
without changing the corresponding disjointness graph.

It is more convenient to speak about intersection graphs rather than their
complement in formulating the next theorem.

{\bf Theorem 13.} {\em
(i)  If $G$ is a graph with maximum degree at most $3$, then $G$ is an intersection graph of lines in $\PP^3$.

(ii) For an arbitrary graph $G$ the {\em line graph} of $G$ is an intersection graph of lines in $\PP^3$.}
\smallskip

\noindent{\bf Proof.} (i)  Suppose first that $G$ is triangle-free. Let
$V(G)=\{v_1,\dots,v_k\}$. Let vertex $v_1$ be represented by an arbitrary line
$\ell_1$. Suppose, recursively, that the line $\ell_j$ representing vertex $j$
has already been defined for every $j<i$. We will maintain the ``general
position'' property that no doubly ruled surface contains more than $3$
pairwise disjoint lines. For the definition and basic properties of doubly ruled surfaces refer to the classic texbook of Hilbert and Cohn-Vossen \cite{HCV}. (The same book can also be used for reference to the few elementary concepts of algebraic geometry we will use below.) We must choose $\ell_i$ representing $v_i$ such that

(a) it intersects the lines representing the neighbors $v_j$ of $v_i$ with $j<i$,

(b) it does not intersect the lines representing the non-neighbors $v_j$ with $j<i$, and

(c) we maintain our general position conditions.

These are simple algebraic conditions. The vertex $v_i$ has at most
$3$ neighbors among $v_j$ for $j<i$, and they must be represented by
pairwise disjoint lines.
%The corresponding
Thus, the Zariski-closed conditions from (a)
determine an irreducible variety of lines, so unless they force the
violation of a specific other (Zariski-open) condition from (b) or (c),
all of those
conditions can be satisfied with a generic line through the lines
representing the neighbors. In case $v_i$ has three neighbors $v_j$
with $j<i$, the corresponding condition forces $\ell_i$ to be in one
of the two families of lines on a doubly ruled surface $\Sigma$. This further
forces $\ell_i$ to intersect {\em all} lines of the other family on $\Sigma$,
but due to the general position condition, none of the vertices of $G$ is
represented by lines there, except the three neighbors of $v_i$. We would
violate the general position condition with the new line $\ell_i$ if the
family we choose it from already had three members representing
vertices. However, this would mean that the degrees of the neighbors of $v_i$
would be at least $4$, a contradiction. In case $v_i$ has fewer than $3$
neighbors, the requirement of $\ell_i$ intersecting the corresponding lines
does not force $\ell_i$ to intersect any further lines or to lie on any doubly
ruled surface.

\smallskip

We prove the general case by induction on $|V(G)|$. Suppose that $a,b,c\in
V(G)$ form a triangle in $G$ and that the subgraph of $G$ induced by
$V(G)\setminus\{a,b,c\}$ can be represented as the intersection graph of
distinct lines in $\PP^3$. Note that each of $a$, $b$ and $c$ has at most a
single neighbor in the rest of the graph. We extend the representation of the
subgraph by adding three lines $\ell_a$, $\ell_b$ and $\ell_c$, representing
the vertices of the triangle. We choose these lines in a generic way so that
they pass through a common point $p$, and $\ell_a$ intersects the line
representing the neighbor of $a$ (in case it exists), and similarly for
$\ell_b$ and $\ell_c$. It is clear that we have enough degrees of freedom (at
least six) to avoid creating any further intersection. For instance, it
suffices to choose $p$ outside all lines in the construction and all planes
determined by intersecting pairs of lines.

\smallskip

(ii) Assign distinct points of $\PP^3$ to the vertices of $G$ so that no four
points lie in a plane. Represent each edge $xx'\in E(G)$ by the line
connecting the points assigned to $x$ and $x'$. As no four points are
coplanar, two lines representing a pair of edges will cross if and only if the
edges share an endpoint. Therefore, the intersection graph of these lines is
isomorphic to the edge graph of $G$. $\Box$

\smallskip

The following theorem implies Theorem~3, as the disjointness graph $H=\bar{G}$
is the complement of the intersection graph $G$, and we have
$\omega(G)=\alpha(H)$, $\alpha(G)=\omega(H)$, $\chi(G)=\theta(H)$, and
$\theta(G)=\chi(H)$. Here $\theta(H)$ denotes the {\em clique covering number}
of $H$, that is, the smallest number of complete subgraphs of $H$ whose vertex
sets cover $V(H)$.
\smallskip

{\bf Theorem 14.} {\em Let $H$ be an intersection graph of
$n$ lines in the Euclidean space $\RR^3$ or in the projective space $\PP^3$.

(i) Computing $\alpha(H)$, the independence number of $H$, is NP-hard.

(ii) Computing $\theta(H)$, the clique covering number of $H$, is NP-hard.

(iii) Deciding whether $\chi(H)\le3$, that is, whether $H$ is $3$-colorable, is NP-complete.

(iv) Computing $\omega(H)$, the clique number of $H$, is in P.

(v) Deciding whether $\theta(H)\le k$ for a fixed $k$ is in P.

(vi) All the above statements remain true if $H$ is not given as an abstract
graph, but with its intersection representation
with lines.}

\medskip

\noindent{\bf Proof.} We only deal with the case where the lines are in
$\PP^3$. The reduction of the Euclidean case to this
case is easy.
\medskip

(i) The problem of determining the independence number of 3-regular graphs is
$NP$-hard; see \cite{AK00}. By Theorem~13(i), all 3-regular graphs are
intersection graphs of lines in $\PP^3$.

\medskip

(ii) The {\em vertex cover number} of a graph $H$ is the smallest number of
vertices with the property that every edge of $H$ is incident to at least one
of them. In
\cite{P74}, it was shown that the problem of determining the {\em
  independendence number} $\alpha(H)$ is $NP$-hard even for triangle-free graphs.
Note that the vertex cover number of $H$ is $|V(H)|-\alpha(H)$.
We can reduce this problem
to the problem of determining the clique covering number of an intersection
graph of lines. For this, note that each complete subgraph of the line graph
$H'$ of a triangle-free graph $H$ corresponds to a star of $H$ and thus $\theta(H')$ is the vertex
cover number of $H$. The reduction is complete, as $H'$ is the intersection
graph of lines in $\PP^3$, by Theorem~13(ii).

\medskip

(iii) Deciding whether the {\em chromatic index} (chromatic number of the line graph)
of a $3$-regular graph is $3$ is
NP-complete, see \cite{H81}. Using that the line graph of any graph
is an intersection graph of lines in $\PP^3$ (Theorem~13(ii)), the statement
follows.

\medskip

(iv) A maximal complete subgraph corresponds to a set of lines passing through
the same point $p$ or lying in the same plane $\Pi$. Any such point $p$ or
plane $\Pi$ is determined by two intersecting lines, so each edge of $H$ is contained in at most $2$ maximal cliques. This limits the number of maximal cliques and yields a polynomial time algorithm to find them all. Computing $\omega(H)$ then reduces simply to finding the largest clique on this list.
\medskip

(v) As we have seen in part (iv), there are at most $2|E(H)|$ we can efficiently find all maximal cliques in $H$. Then we can check all $k$-tuples of them to decide whether they cover all vertices in $H$.

\medskip

(vi) For the statements (iv--v) the claim is trivial: we can ignore the extra information given in the input.

To see the claim for the statements (i--iii), we need to consider the constructions of lines in the
representations described in the proof of Theorem~13, and show that they can be
built in
polynomial time. This is obvious in part (ii) of the theorem. For part (i),
the situation is somewhat more complex. To find many possible representations
of the next vertex intersecting the lines it should, is an algebraically simple task.
In polynomial time, we can find one of them that is generic in the sense needed
for the construction. However, if the coordinates of each line would be twice
as long as those of the preceding line (a condition that is hard to rule out {\em a
priori}), then the whole construction takes more than polynomial time.

A simple way to avoid this problem is the following. First, color the vertices
of the triangle-free graph $G$ of maximal degree at most $3$ by at most $4$
colors, by a simple greedy algorithm.
%
%(or even by $3$) colors.
Find the lines representing the vertices in the following
order: first for the first color class, next for second color class, etc. The
coordinates of each line will be just slightly more complex than the
coordinates of the lines representing vertices in {\em earlier color classes}. Therefore, the construction can be performed in polynomial time. A similar
argument works also for graphs $G$ with triangles: First we find a maximal
subset of pairwise vertex-disjoint triangles in $G$. Let $G_0$ be the graph
obtained from $G$ by removing these triangles. Then we construct an auxiliary
graph $G'$ with these triangles as vertices by connecting two of them with an edge
if there is an edge in $G$ between the triangles. The graph $G'$ has maximum
degree at most $3$, so it can be greedily $4$-colored. If we construct $G$ by
adding back the triangles to $G_0$, in the order determined
by their colors, then the procedure will end in polynomial time.
$\Box$
\medskip

\noindent{\bf Acknowledgement.} A preliminary version of this paper with the same title was published in the {\em Proceedings of the 33rd International Symposium on Computational Geometry} 2017: 59:1--59:15.

\end{document}